\documentclass[11pt]{amsart}
 \usepackage{amsthm,amsmath,amssymb,amscd,epsfig}

 {\end{list}}
 {
    \newtheorem{theorem}{Theorem}[section]
    \newtheorem{proposition}[theorem]{Proposition}
    \newtheorem{lemma}[theorem]{Lemma}

    \newtheorem{corollary}[theorem]{Corollary}

    \newtheorem{remark}[theorem]{Remark}
 }
 {\theoremstyle{definition}

    \newtheorem{definition}[theorem]{Definition}
 }

 \newcommand{\av}{{\bf a}}
 \newcommand{\bv}{{\bf b}}
 \newcommand{\cv}{{\bf c}}
 \newcommand{\dv}{{\bf d}}
 \newcommand{\ab}{\av\bv}
 \newcommand{\cd}{\cv\dv}
 \newcommand{\hz}{{\widehat{0}}}
 \newcommand{\ho}{{\widehat{1}}}
 \newcommand{\join}{\vee}
 \newcommand{\meet}{\wedge}
 
 \newcommand{\zcd}{\ZZ\langle \cv,\dv \rangle}

 \newcommand{\onethingatopanother}[2]
 {\genfrac{}{}{0pt}{}{#1}{#2}}

 \newcommand{\RR}{{\mathbb{R}}}

 \newcommand{\ZZ}{{\mathbb{Z}}}

 \newcommand{\cC}{{\mathcal C}}
 \newcommand{\cD}{{\mathcal D}}
 \newcommand{\cF}{{\mathcal F}}

 \newcommand{\cL}{{\mathcal L}}

 \newcommand{\cO}{{\mathcal O}}

 \newcommand{\Int}{\operatorname{Int}}
 \newcommand{\Pyr}{\operatorname{Pyr}}
 \newcommand{\wt}{\operatorname{wt}}

 \newcommand{\res}{\operatorname{res}}
 
 \newcommand{\Link}{\operatorname{Link}}

 \newcommand{\lar}{\longleftarrow}
 \newcommand{\rar}{\longrightarrow}
 
 \newcommand{\isom}{\simeq}

 \newcommand{\Cd}{C^\bullet}
 \newcommand{\Po}{Poincar\'e }
 \newcommand{\CM}{Cohen-Macaulay}
 \newcommand{\setmin}{{\smallsetminus}}

\begin{document}
 \title{Decomposition theorem for the $\cd$-index of Gorenstein* posets}

 \author{Richard Ehrenborg}
 \address{Department of Mathematics\\
 University of Kentucky\\
 Lexington, KY 40506\\
 USA.}
 \email{jrge@ms.uky.edu}

 \author{Kalle Karu}
 \address{Department of Mathematics\\ University of British Columbia \\
   1984 Mathematics Road\\
 Vancouver, B.C. Canada V6T 1Z2}
 \email{karu@math.ubc.ca}

 \begin{abstract}
 We prove a decomposition theorem for the $\cd$-index of a
 Gorenstein* poset analogous to the decomposition theorem for the
 intersection cohomology of a toric variety. From this we settle a
 conjecture of Stanley that the $\cd$-index of Gorenstein*
 lattices is minimized on Boolean algebras.
 \end{abstract}

 \maketitle

 \section{Introduction}
 \setcounter{equation}{0}

 The $\cd$-index of a convex polytope $P$ is a polynomial $\Psi_{P}(\cv,\dv)$
 in the non-commuting variables
 $\cv$ and $\dv$ that effectively encodes the flag $f$-vector of
 the polytope $P$~\cite{BB,BK}.
 Its coefficients are non-negative integers~\cite{S}.

 The following result was proved in \cite{BE} and used to study the
 monotonicity property of the $\cd$-index:

 \begin{theorem}[Billera, Ehrenborg] \label{thm-BE}
 For any polytope $P$ and a proper face $F$ of the polytope $P$
 the coefficientwise inequality
 \[ \Psi_P \geq \Psi_F \cdot \Psi_{\Pyr(P/F)}\]
 holds, where $\Pyr(P/F)$ is the pyramid over the polytope $P/F$.
 \label{theorem_Billera_Ehrenborg}
 \end{theorem} 
 By iterating this theorem it implies that among all polytopes of
 dimension $n$ the $\cd$-index is minimized by the $n$-dimensional
 simplex.

 One can define the $\cd$-index more generally for Eulerian posets. It
 has non-negative coefficients if the poset is
 Gorenstein*~\cite{K}. We generalize
 Theorem~\ref{theorem_Billera_Ehrenborg}
 to the case of Gorenstein* lattices:

 \begin{theorem} \label{thm-main} 
 For any Gorenstein* lattice $\Lambda$ and an element
 $\nu \in \Lambda$ such that $\hz < \nu < \ho$ 
 the coefficientwise inequality
 \[ \Psi_\Lambda \geq \Psi_{[\hz,\nu)} \cdot \Psi_{\Pyr[\nu,\ho)}\]
 holds, where $\Pyr[\nu,\ho)$ is the pyramid over the poset
   $[\nu,\ho)$.
 \end{theorem} 

 The use of half-open intervals in the statement of the theorem is
 explained in the next section. 
 Considering the dual lattice
 $\Lambda^{*}$, we obtain by duality
 \[  \Psi_\Lambda \geq \Psi_{\Pyr[\hz,\nu)} \cdot
   \Psi_{[\nu,\ho)}.\]
 The action of taking the pyramid
 of an Eulerian poset on the $\cd$-index is described
 by the following linear map.
 Let $G$ be a derivation on $\zcd$
 defined by $G(\cv) = \dv$ and $G(\dv) = \cv\dv$.
 Next let $\Pyr$ be the operator on $\zcd$
 defined by $\Pyr(w) = w \cdot \cv  + G(w)$.
 Then the $\cd$-index of the pyramid is the pyramid of the $\cd$-index,
 that is, $\Psi_{\Pyr(P)} = \Pyr(\Psi_{P})$;
 see~\cite{ER}.

 As an example, if $\Lambda$ is the face lattice of a plane $n$-gon
 and $\nu$ is a vertex, then the inequality
 of Theorem~\ref{thm-main} reads
 \[ \cv^2+(n-2)\dv \geq \cv^2+\dv.\]
 To see that the theorem does not extend to Gorenstein* posets that are
 not lattices, take $\Lambda$ to be the poset of the plane
 $2$-gon, that is, the $CW$-complex consisting of two edges glued at the
 two endpoints. The cells of this complex form a Gorenstein* poset, but the
 inequality above with $n=2$ does not hold.

 As a corollary to Theorem~\ref{thm-main} we settle
 a conjecture of Stanley:

 \begin{corollary} Among all Gorenstein* lattices the $\cd$-index is
   minimized by the Boolean algebra.
 \end{corollary}

 Recall that to a polytope $P$ one can associate its $h$-polynomial
 $h_P(t)$ and its $g$-polynomial $g_P(t)$. The inequality in
 Theorem~\ref{thm-BE} was motivated by a similar inequality between the
 $g$-polynomials conjectured by Kalai and proved in the case of
 rational polytopes by Braden and
 MacPherson \cite{BM}.

 \begin{theorem}[Braden, MacPherson] \label{thm-BM}
  For any rational polytope $P$ and a proper face $F$ of $P$
 the following inequality holds coefficientwise:
 \[ g_P \geq g_F \cdot g_{P/F}  . \]
 \end{theorem} 

 Since $g_{P/F} = g_{\Pyr(P/F)}$, the two inequalities in
 Theorems~\ref{thm-BE} and \ref{thm-BM} can be made to
 look even more similar. 

 Theorem~\ref{thm-BM} was generalized to the case of nonrational
 polytopes in \cite{BBFK, BL}, subject to the assumption that the
 $g$-polynomials involved are non-negative \cite{K0}. This
 generalization is proved using a combinatorial decomposition theorem
 for pure sheaves on a fan: a pure (i.e., locally free and flabby)
 sheaf on a fan decomposes into elementary sheaves. Another consequence
 of this combinatorial decomposition theorem is the monotonicity of the
 $h$-vector under subdivisions: if $\widehat{\Delta}$ is a subdivision
 of a complete fan $\Delta$ then
 \[ h_{\widehat{\Delta}} \geq h_\Delta\]
 coefficientwise.

 We prove an analogous decomposition theorem for the $\cd$-index. A more
 precise statement is given in Section~2.

 \begin{theorem} \label{thm-decomp}
 Let $\widehat{\Pi}$ be a subdivision of a Gorenstein* poset $\Pi$. Then 
 the following inequality holds coefficientwise:
 \[ \Psi_{\widehat{\Pi}} \geq \Psi_{\Pi}   .  \]
 \end{theorem}

 The crucial point in the last theorem is the correct definition of
 ``subdivision''. It not only
 includes the usual subdivisions of polyhedral
 fans, but also more general subdivisions of $CW$-complexes. With a
 correct notion of subdivision in place, Theorem~\ref{thm-main} follows
 from Theorem~\ref{thm-decomp}.  The product of $\cd$-indices
 corresponds to the $*$-product of posets (see
 Section~\ref{subsec-lattices} below for the definition of $*$-product):
 \[ \Psi_{[\hz,x)} \cdot \Psi_{\Pyr[x,\ho)} =
     \Psi_{[\hz,x)*\Pyr[x,\ho)}.\]
 We show that the original lattice $\Lambda$ is a subdivision of
 $[\hz,x)*\Pyr[x,\ho)$ and then apply Theorem~\ref{thm-decomp}.

\begin{figure}[ht] \label{fig-pol}
 \centerline{\psfig{figure=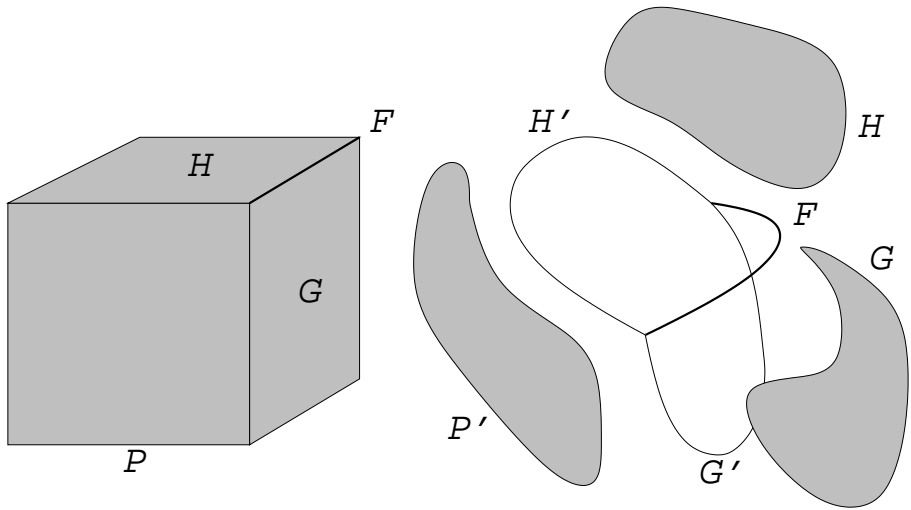,width=12cm}}
 \caption{Construction of $F * \Pyr(P/F)$ from $P$ and $F$.}
 \end{figure}

 For a polytope $P$, the construction of $F * \Pyr(P/F)$ is illustrated in
 Figure~\ref{fig-pol}. From the original polytope $P$ we keep all faces
 $G\leq F$ and $G> F$. Then for each $G>F$ we cap it off with a cell
 $G'$ of one smaller dimension. The resulting $CW$-complex corresponds
 to the poset $F * \Pyr(P/F)$ and it admits (the boundary complex of)
 $P$ as its subdivision.

 \section{Decomposition theorem} 
 \setcounter{equation}{0}

 We refer to \cite{EC1,S1} for terminology about posets.
 Throughout we only consider
 posets that are finite, graded, with minimal element $\hz$, but with
 no maximal element $\ho$ in general. This is motivated by the
 poset of a fan: every fan has a minimal cone $\{0\} = \hz$, but no maximal
 cone in general. The rank of an element $x$ is
 denoted by $\rho(x)$. If $x<y$, then $\rho(x,y) = \rho(y)-\rho(x)$. We
 assume that the minimal element has rank 
 $\rho(\hz)=0$ and maximal elements have rank equal to the rank of
 the poset.

 The closed  (resp. half-open) intervals will be denoted $[x,y]$
 (resp. $[x,y)$). If there is possibility of confusion, we write the
 poset $\Pi$ as the subscript: $[x,y]_\Pi$. Even though $\Pi$ may
 not contain  the maximal element $\ho$, we write $[x,\ho)$ for
 all the elements greater than or equal to $x$.

 If $\Pi$ is a poset of rank $n$, let $\cO(\Pi)$ be the order complex,
 that is,
 the set of chains in $\Pi$ ordered by inclusion.
 \[  \cO(\Pi) = \{ \{\hz=\sigma_{0}<\sigma_1<\sigma_2<\cdots<\sigma_k\} |
 \sigma_i\in\Pi, k \geq 0\}.\]
 Then  $\cO(\Pi)$ is again a poset of rank $n$. In fact, $\cO(\Pi)$ is
 the face poset of a simplicial complex of dimension $n-1$. (More
 precisely, the subsets  of $\Pi$ of the form
 $\{\sigma_1,\sigma_2,\ldots,\sigma_k\}$, where
 $\hz<\sigma_1<\sigma_2<\cdots<\sigma_k$ form a simplicial complex.)

 If a poset $\Pi$ is a lattice, then it obviously must have a maximal
 element $\ho$. If $\Pi$ does not contain $\ho$, then by abuse of
 terminology we say that $\Pi$ is a lattice of rank $n$ if
 $\Pi\cup\{\ho\}$ is a lattice of rank $n+1$. For instance, the
 simplicial complex  $\cO(\Pi)$ is a lattice of rank $n$ for any poset
 $\Pi$ of rank $n$.

 A poset $\Pi$ is Eulerian if every interval satisfies the
 Euler-Poincar\'e relation. Equivalently, for any two elements
 $\tau$, $\pi$ in the poset $\Pi$ such that
 $\hz \leq \tau < \pi \leq \ho$, 
 we have 
 $$   \sum_{\tau \leq \sigma \leq \pi} (-1)^{\rho(\tau,\sigma)} = 0  .  $$

 \subsection{Gorenstein* posets}

 \begin{definition} \label{def-Gor} A simplicial complex $\Delta$ of
   dimension  $n-1$ is called {\em Gorenstein} if it is a real homology
   sphere of dimension $n-1$. This means that:
 \begin{itemize}
 \item The reduced (simplicial) homology of $\Delta$ is 
 \[ \widetilde{H}_i(\Delta,\RR) = \begin{cases} \RR & \text{if $i=n-1$}\\
 0 & \text{otherwise.} \end{cases}\]
 \item For every simplex $x\in\Delta$ of dimension $m$, the
   reduced homology of its link is
  \[ \widetilde{H}_i(\Link_\Delta(x),\RR) = \begin{cases} \RR & \text{if $i=n-m-1$}\\
 0 & \text{otherwise.} \end{cases}\]
 \end{itemize}
 A poset $\Pi$ is called {\em Gorenstein*} if its order complex
 $\cO(\Pi)$ is Gorenstein. 
 \end{definition}

 A simplicial complex $\Delta$ is called {\CM} if a  similar condition
 holds, but the top degree cohomology is allowed to have any
 dimension. Thus a simplicial
 complex is Gorenstein if and only if it is {\CM} and Eulerian. A poset
 $\Pi$ is called {\CM} if its order complex $\cO(\Pi)$ is {\CM}. Since
 $\Pi$ is Eulerian if and only if $\cO(\Pi)$ is Eulerian, it follows
 that a poset $\Pi$ is Gorenstein* if and only if it is {\CM} and
 Eulerian. 

 Recall that $\Link_\Delta(x)$ is isomorphic to the interval $[x,\ho)$
 in $\Delta$. Thus, denoting $\Link_\Delta(\hz)=\Delta$, we see that
 the first homology condition in the definition is the same as the
 second one applied to the element $x=\hz$. Also recall that the
 reduced simplicial homology of $\Link_\Delta(x)\isom [x,\ho)$ for $x$
   of rank $k$ is computed by the simplicial chain complex
 \[   0   \lar
    \RR_x \lar
    \bigoplus_{\onethingatopanother{y>x}{\rho(y)=k+1}} \RR_y \lar
    \bigoplus_{\onethingatopanother{z>x}{\rho(z)=k+2}} \RR_z \lar 
    \cdots \lar
    \bigoplus_{\onethingatopanother{w>x}{\rho(w)=n}} \RR_w \lar 0  .
 \]
 The homology condition in the definition of Gorenstein states that
 this sequence has only $1$-dimensional homology at the rightmost place.

 If $\Pi$ is a Gorenstein* poset of rank $n$, then any interval $[\sigma,\tau)$
   in $\Pi$ is again Gorenstein*. Indeed, given a chain $x\in
   \cO([\sigma,\tau))$, the interval $[x,\ho)$ in
   $\cO([\sigma,\tau))$ is isomorphic to the interval $[y,\ho)$ in
   $\cO(\Pi)$ for some $y\in \cO(\Pi)$. Hence the homology condition
   for $[y,\ho)$ implies the same for $[x,\ho)$.

 Given a Gorenstein* poset $\Pi$, its dual poset $\Pi^{*}$ (obtained
 by reversing the order relations) is again Gorenstein*. This follows
 from the fact that the two posets have isomorphic order complexes.

 We also need the notion of a poset $\Pi$ such that $\cO(\Pi)$ is a
 homology ball. If such an $\cO(\Pi)$ has dimension $n-1$, then its
 boundary $\partial\cO(\Pi)$ is a sub-complex of dimension $n-2$. To
 define the analogue for arbitrary posets, we consider pairs
 $(\Pi,\partial\Pi)$, where
 \begin{itemize}
 \item $\Pi$ is a poset of rank $n$;
 \item $\partial\Pi\subset \Pi$ is a sub-poset of $\Pi$ of rank $n-1$
   (with elements in ranks $n-1$ or less);
 \item $\partial\Pi$ is an ideal in $\Pi$: if $\sigma\in \partial\Pi$
   and $\tau<\sigma$, then $\tau\in \partial\Pi$.
 \end{itemize}

 \begin{definition} \label{def-nGor}
 A pair $(\Delta,\partial\Delta)$ where $\Delta$ is a simplicial
 complex of dimension $n-1$ is called {\em near-Gorenstein} if $\Delta$
 is a real homology ball of dimension $n-1$ with
 boundary~$\partial\Delta$. This means that:
 \begin{itemize}
 \item The complex $\partial\Delta$ is Gorenstein of dimension $n-2$.
 \item For every $x\in\Delta$ of dimension $m$ (including $x=\hz$ of
   dimension $-1$), the reduced homology of its link is
  \[ \widetilde{H}_i(\Link_\Delta(x),\RR) = \begin{cases} \RR & \text{if
   $i=n-m-1$ and $x\notin \partial\Delta$}\\
 0 & \text{otherwise.} \end{cases}\]
 \end{itemize}
 A pair $(\Pi,\partial\Pi)$ is called {\em near-Gorenstein*} if 
 $(\cO(\Pi), \cO(\partial\Pi))$ is near-Gorenstein.
 \end{definition}

 It is clear that for a near-Gorenstein* pair $(\Pi,\partial\Pi)$, the
 boundary $\partial\Pi$ is determined by $\Pi$. Hence we may simply say
 that a poset $\Pi$ is near-Gorenstein*. Similarly, we may call a
 simplicial complex $\Delta$ near-Gorenstein. We also denote $\Int(\Pi) =
 \Pi\setmin\partial\Pi$.

 The terminology is motivated by Stanley's notion of a near-Eulerian
 poset \cite{S}, which is a poset obtained from an Eulerian poset of
 rank $n$ by removing an element of rank $n$. A similar result holds here:

 \begin{lemma} \label{lem-Gor-near-Gor}
 Let $\Pi$ be a Gorenstein* poset of rank $n$ and $\pi\in \Pi$ an
 element of rank $\rho(\pi) = n$. Then $\Pi\setmin\{\pi\}$ is a
 near-Gorenstein* poset of rank $n$ with boundary $[\hz,\pi)$.

  Conversely, every near-Gorenstein* poset arises this way.
 \end{lemma}

 We will prove this lemma in Section~\ref{section_sheaves}.

 \subsection{The cd-index}\label{sec-cd-index}

 Let us recall the combinatorial construction of the $\cd$-index of an
 Eulerian poset. For a chain
 $x = \{\hz = \sigma_{0} < \sigma_{1} < \cdots < \sigma_{k}\}$
 in a poset $\Pi$ define the weight of the chain
 by
 $$   \wt(x)
   =
      (\av - \bv)^{\rho(\sigma_{0},\sigma_{1})-1}
     \cdot
      \bv
     \cdot
      (\av - \bv)^{\rho(\sigma_{1},\sigma_{2})-1}
     \cdot
      \bv
     \cdots
      \bv
     \cdot
      (\av - \bv)^{\rho(\sigma_{k},\ho)-1}  .   $$

 Here $\av$ and $\bv$ are non-commuting variables of degree $1$ each.
 Then the $\ab$-index of the poset $\Pi$ is given by the
 sum
 $$    \Psi_{\Pi} = \Psi_{\Pi}(\av,\bv) = \sum_{x} \wt(x)  ,  $$
 where the sum ranges over all chains $x$ in the poset $\Pi$.


 When $\Pi$ is Eulerian, it follows from the relations in  \cite{BB}
 that $\Psi_\Pi(\av,\bv)$ can be expressed in terms of the variables
 $\cv=\av+\bv$
 and $\dv = \av\bv + \bv\av$. The $\cd$-index of~$\Pi$ is the polynomial
 $\Psi_\Pi(\av,\bv)$ written in terms of $\cv$ and $\dv$. By a slight
 abuse of notation we denote the $\cd$-index of $\Pi$ by
 $\Psi_\Pi(\cv,\dv)$. This notation may cause concern only when we talk
 about non-negativity of the  $\cd$-index. By this we mean
 non-negativity of the coefficient of $\Psi_\Pi(\cv,\dv)$ when
 expressed as a polynomial in  $\cv$ and $\dv$. (The
 coefficients of  $\Psi_\Pi(\av,\bv)$ are always non-negative.)

  Note that
 $\Psi_\Pi(\av,\bv)$ is homogeneous of degree~$n$. The same is true for
 $\Psi_\Pi(\cv,\dv)$ where $\cv$ has degree $1$ and $\dv$ has degree $2$.

 The following non-negativity theorem proved in~\cite{K}
 by the second author generalizes
 Stanley's result of the non-negativity for the $\cd$-index of
 $S$-shellable complexes~\cite{S}.

 \begin{theorem} \label{thm-pos-cd} The $\cd$-index of a Gorenstein*
   poset has non-negative integer coefficients.
 \end{theorem}

 One can also define the $\ab$-index of a near-Gorenstein* poset
 $\Pi$ using the same construction. This $\ab$-polynomial, however,
 cannot be expressed in terms of $\cv$ and~$\dv$. Considering chains in
 the Gorenstein* poset $\Pi\cup\{\pi\}$ (Lemma~\ref{lem-Gor-near-Gor}),
 it is easy to see that:
 \[ \Psi_\Pi(\av,\bv) = \Phi_\Pi(\av,\bv) +\Psi_{\partial\Pi}(\av,\bv)\cdot \av,\]
 where $\Phi_\Pi(\av,\bv)$ is a homogeneous $\ab$-polynomial that can be
 expressed in terms of $\cv$ and $\dv$. Let us omit the multiplication with
 $\av$ at the end and define the $\cd$-index of $\Pi$ to be
 \[ \Psi_\Pi(\cv,\dv) := \Phi_\Pi(\cv,\dv)
 +\Psi_{\partial\Pi}(\cv,\dv), \label{eq-cd-nGor}\]
 where $\Phi_\Pi(\cv,\dv)$ is the polynomial $\Phi_\Pi(\av,\bv)$
 expressed in terms of  $\cv$ and $\dv$.
 Note that the first summand is homogeneous of degree $n$ and the
 second summand is homogeneous of degree $n-1$.

 The following Theorem is stated in \cite{K} only for the
 case of fans. We give a generalization to the case of
 near-Gorenstein* posets in the appendix.

 \begin{theorem} \label{thm-pos-cd-nGor} The $\cd$-index of a
   near-Gorenstein* poset has non-negative integer coefficients.
 \end{theorem}

 \subsection{Subdivisions}

 Let us now define the notion of a subdivision.

 \begin{definition} \label{def-subd}
 Let $\widehat{\Pi}$ and $\Pi$ be two Gorenstein* posets of rank $n$ and 
 let
 \[ \phi: \widehat{\Pi}\to \Pi \]
 be a surjective map that preserves the order relation, but not
 necessarily the rank. Then $\phi$ is a {\em subdivision} if for
 every $\sigma\in\Pi$ the pair
 \[ (\phi^{-1} [\hz,\sigma], \phi^{-1} [\hz,\sigma))\]
 is near-Gorenstein* of rank equal to $\rho(\sigma)$.
 \end{definition}

 To simplify notation, let us write $\widehat{\sigma} = \phi^{-1}
 [\hz,\sigma]$ and $\partial\widehat{\sigma} = \phi^{-1}
 [\hz,\sigma)$. Recall that by Theorem~\ref{thm-pos-cd-nGor}, the
   $\cd$-index 
 \begin{equation} \Psi_{\widehat{\sigma}} = \Phi_{\widehat{\sigma}}
   +\Psi_{\partial\widehat{\sigma}} \label{eq-sigma} \end{equation}
 has all non-negative coefficients. 

 We can now state the decomposition theorem.

 \begin{theorem} \label{thm-decomp-prec}
 Let $\phi:\widehat{\Pi}\to \Pi$ be a subdivision map. Then
 \[ \Psi_{\widehat{\Pi}} = \sum_{\sigma\in\Pi}
 \Phi_{\widehat{\sigma}}\cdot\Psi_{[\sigma,\ho)}.\]
 \end{theorem}

 \begin{proof} The proof of this theorem is almost entirely contained
   in \cite{BBFK, BL, K}. It relies on the theory of pure sheaves
  on fans and posets. Since we will
   not use this theory elsewhere in this article, we will use the
   notation of the references above and sketch the proof.

 The category of pure sheaves (i.e., locally free and flabby sheaves)
 on a fan $\Delta$ is semisimple: any pure sheaf $\cF$ has a
 decomposition into a direct sum of simple sheaves $\cL_\sigma$ indexed
 by cones $\sigma\in\Delta$:
 \[ \cF = \bigoplus_{\sigma\in\Delta} V_\sigma\otimes \cL_\sigma.\]
 Here $V_\sigma$ is a graded vector space,
 \[ V_\sigma \isom \ker(\overline{\cF_\sigma} \rar
 \overline{\cF_{\partial\sigma}}).\]
 Since the map $\overline{\cF_\sigma} \rar
 \overline{\cF_{\partial\sigma}}$ is surjective by purity of $\cF$, we
 get an equality between the \Po polynomials of the graded vector
 spaces:
 \[ P_{V_\sigma}(t) = P_{\overline{\cF_\sigma}}(t) -
 P_{\overline{\cF_{\partial\sigma}}}(t).\]

 It is shown in \cite{K} that one can define a similar theory of pure
 sheaves on an arbitrary poset $\Pi$, where the sheaves now are
 multi-graded by $\ZZ^n$. If $\Pi$ is Gorenstein* or near-Gorenstein* then
 the $\cd$-index of $\Pi$ is obtained by a change of variable formula
 from the multi-variable \Po polynomial of the
 module of global sections of $\cL_\hz$. For $\sigma>\hz$, the sheaf
 $\cL_\sigma$ similarly gives the $\cd$-index of $[\sigma,\ho)$.

 If $\phi:\widehat{\Pi} \to \Pi$ is a subdivision map, we can decompose 
 \[ \phi_*(\cL_{\hz}) = \bigoplus_{\sigma\in\Pi} V_\sigma \otimes
 \cL_\sigma.\] 
 Here 
 \[ V_\sigma \isom \ker(\overline{\cL_{\widehat\sigma}} \rar
 \overline{\cL_{\partial{\widehat\sigma}}}).\]
 Its multi-variable \Po polynomial is
 \[  P_{V_\sigma}({\bf t}) = P_{\overline{\cL_{\widehat\sigma}}}({\bf t}) -
 P_{\overline{\cL_{\partial{\widehat\sigma}}}}({\bf t}) =
   (\Phi_{\widehat{\sigma}} +\Psi_{\partial\widehat{\sigma}}) -
   \Psi_{\partial\widehat{\sigma}} = \Phi_{\widehat{\sigma}},\]
 where the middle equality comes from Equation~(\ref{eq-sigma}).
 Now the formula in the theorem follows by taking
 the \Po polynomials of the modules of global sections on both sides of
 the decomposition.
 \end{proof}

 Note that all of the terms in the sum of Theorem~\ref{thm-decomp-prec}
 are non-negative. Since $\Phi_{\hz}=1$, we obtain
 Theorem~\ref{thm-decomp} as a corollary.

 \subsection{Lattices} \label{subsec-lattices}

 Let us now deduce Theorem~\ref{thm-main} from
 Theorem~\ref{thm-decomp}. To do this, we construct a subdivision map
 $\phi: \Lambda\to [\hz,\nu)*\Pyr[\nu,\ho)$.

 Given two posets $\Pi$ and $\Lambda$, define 
 \[ \Pi * \Lambda = \Pi\cup (\Lambda\setmin\{\hz\}),\]
 with order relations those of $\Pi$ and $\Lambda$,  and additionally
 $\sigma<\tau$ for any $\sigma\in\Pi$ and $\tau\in\Lambda$. If $\Pi$ and
 $\Lambda$ have ranks $m$ and $n$, respectively, then $\Pi * \Lambda$
 has rank $m+n$. The order poset $\cO(\Pi * \Lambda)$ is the product
 $\cO(\Pi) \times \cO(\Lambda)$. It follows that if $\Pi$ and $\Lambda$
 are both Gorenstein*, then  so is $\Pi * \Lambda$. Also, if $\Pi$
 is Gorenstein* and $\Lambda$ is near-Gorenstein*, then $\Pi *
 \Lambda$ is near-Gorenstein*. These two properties follow by
 considering tensor products of the chain complexes of $\cO(\Pi)$ and
 $\cO(\Lambda)$. Moreover, the $\cd$-index turns
 the $*$-product into multiplication~\cite{S}:
 \[ \Psi_{\Pi * \Lambda} = \Psi_{\Pi} \cdot \Psi_{\Lambda}.\]

 The pyramid operation of a poset containing the
 maximal element $\ho$ is defined as the
 Cartesian product with the Boolean algebra $B_{1} = \{\hz,\ho\}$.
 That is, for our poset $\Pi$
 we have that
 $\Pyr(\Pi) \cup \{\ho\} = (\Pi\cup \{\ho\})\times B_{1}$.
 Let $P$ and $Q$ be two posets both containing maximal elements.
 Directly we have that if both $P$ and $Q$ are Eulerian
 so is their Cartesian product $P \times Q$.
 Similarly, if both $P$ and $Q$ are {\CM} posets
 it follows from \cite[Theorem~7.1]{B} that
 their product $P \times Q$ is {\CM}.
 Hence the class of Gorenstein* posets
 are closed under Cartesian product
 and we conclude that:
 \begin{lemma} Let $\Pi$ be a Gorenstein* poset of rank $n$. Then
   $\Pyr(\Pi)$ is a Gorenstein* poset of rank $n+1$.
 \end{lemma}

 With notation as in Theorem~\ref{thm-main}, let $\Lambda$ be a
 Gorenstein* lattice and $\nu\in\Lambda$, $\hz<\nu<\ho$.
 Define the map $\phi: \Lambda\to
 [\hz,\nu)*\Pyr[\nu,\ho)$ as follows:
 \[ \phi(\tau) = \begin{cases} \tau & \text{if $\tau<\nu$}\\
 (\tau,\ho) & \text{if $\tau\geq\nu$}\\
 (\tau \vee \nu,\hz) & \text{otherwise}.
 \end{cases}\]

 Note that $\tau\vee\nu$ can be equal to $\ho$. In
 Figure~\ref{fig-pol} the cells $G$ and $G'$ correspond to elements
 $(G,\ho)$ and $(G,\hz)$, and the cell $P'$ corresponds to
 $(\ho,\hz)$.

 We claim that $\phi$ is a subdivision map. The map $\phi$ clearly
 preserves the order relation. To prove surjectivity of $\phi$, the
 only nontrivial case is to show that $(\sigma,\hz)$ lies in the image
 for $\sigma > \nu$. Replacing $\Lambda$ by the interval $[\hz,
 \sigma)$, we need to show that for some $\tau\in\Lambda$, we have
 $\tau\vee\nu = \ho$. If no such $\tau$ exists, then every maximal
 element of $\Lambda$ lies in $[\nu,\ho)$. By descending induction on
 the rank of an element, and using that the lattice is Eulerian,
 it follows that every element must lie in 
 $[\nu,\ho)$. This gives a contradiction.

 Finally, we  need to check that the
 inverse images $(\widehat{\pi}, \partial\widehat{\pi})$ are
 near-Gorenstein*  in the three cases: 
 \begin{enumerate}
 \item If $\pi=\tau<\nu$, then 
 \[ (\widehat{\pi}, \partial\widehat{\pi}) = ([\hz,\tau],[\hz,\tau)).\]
 Let $\Pi$ be the Gorenstein* poset $[\hz,\tau)$. We need to show
 that $\Pi\cup\{\ho\}$ is near-Gorenstein* with boundary $\Pi$.

 \item If $\pi =(\tau,\ho)$, then we have as in the previous case
 \[  (\widehat{\pi}, \partial\widehat{\pi}) = ([\hz,\tau],[\hz,\tau)).\]

 \item If $\pi =(\tau,\hz)$, then
 \[ \widehat{\pi} = [\hz,\tau)\setmin [\nu,\ho),\]
 with boundary 
 \[ \partial\widehat{\pi} = \{\sigma \in \widehat{\pi}| \sigma\vee \nu < \tau\}.\]
 Let $\Pi$ be the Gorenstein* lattice $[\hz, \tau)$ containing
   $\nu$ (recall that by this we mean $\Pi\cup\{\ho\}$ is a
   lattice). Then we need to show that $\Pi\setmin [\nu,\ho)$ is 
   near-Gorenstein* with boundary 
 \[  \{\sigma| \sigma\vee \nu < \ho \}.\]
 \end{enumerate}

 We state what is left to prove in the following two lemmas.

 \begin{lemma}
 \label{lem-add-one}
 Let $\Pi$ be a Gorenstein* poset of
 rank $n$. Then 
 $\Pi\cup\{\ho\}$ is a near-Gorenstein* poset of rank $n+1$ with
 boundary $\Pi$.
 \end{lemma}

 \begin{proof}
 Note that  $\Pi\cup\{\ho\}$ is the $*$-product 
 $\Pi *\{\hz,\ho\}$,
 where $\{\hz,\ho\}$ is near-Gorenstein* with boundary
 $\{\hz\}$. This implies the statement of the lemma.
 \end{proof}

 \begin{lemma}
 \label{lem-main}
 Let $\Lambda$ be a Gorenstein* lattice
 of rank $n$ and 
 $\nu\in\Lambda\setmin{\hz}$. Then 
 \[ \Lambda\setmin [\nu,\ho) \]
 is near-Gorenstein* of rank $n$ with boundary 
 \[ \{ \tau\in\Lambda\setmin[\nu,\ho) | \tau\vee\nu < \ho \}.\]
 \end{lemma}

 We will prove this lemma in 
 Section~\ref{section_sheaves}.

 As a final remark in this section, observe that we used
 the join operation $\join$ in the lattice, but not the meet
 operation $\meet$. 
 It is enough to do this since
 a finite join-semilattice
 with a minimal element $\hz$ is a lattice;
 see~\cite[Proposition~3.3.1]{EC1}.


 \section{Flag enumeration}
 \setcounter{equation}{0}

 In this section we consider Eulerian posets and lattices
 that are not necessarily Gorenstein*.
 By flag enumeration we are still able to
 derive identities for their $\cd$-indexes.

 \begin{definition}
 For an Eulerian lattice $\Lambda$ with an element $\nu$ such
 that $\hz < \nu < \ho$ define the two subposets
 $\Lambda_{\nu}$ and $\Lambda_{\nu}^{\prime}$
 of the lattice $\Lambda$ as follows. Let
 $$  \Lambda_{\nu}
   = 
     \{ \sigma \in \Lambda | \sigma \join \nu < \ho \} , $$
 where $\Lambda_{\nu}$ inherits the order relations from $\Lambda$,
 and let $\Lambda_{\nu}^{\prime}$ be the
 semisuspension of~$\Lambda_{\nu}$, that is,
 $$  \Lambda_{\nu}^{\prime}
   = 
     \Lambda_{\nu}
   \cup
     \{ * \} ,   $$
 where the new element $*$ satisfies
 the order relation
 $\sigma < *$ for all $\sigma \in \Lambda_{\nu}$ such that $\nu
   \not\leq \sigma$. 
 \end{definition}

 The geometric intuition for the poset $\Lambda_{\nu}$ is that
 it is the near-Eulerian poset consisting of all
 faces contained in facets which contain the face $\nu$. If $\Lambda$
 is Gorenstein*, then we will show in Section~\ref{section_sheaves}
 that $\Lambda_\nu$ is near-Gorenstein*. In fact, the posets $\Lambda\setmin
 [\nu,\ho)$ of Lemma~\ref{lem-main} and $\Lambda_{\nu}$ are
 complementary in that the lattice $\Lambda$ is obtained by gluing the two
 posets along their common boundary.

 \begin{proposition}
 The $\ab$-index of the poset $\Lambda_{\nu}$ is given by
 $$  \Psi_{\Lambda_{\nu}}
   =
     \sum_{\nu \leq \pi < \ho}
		\Psi_{[\hz,\pi)} \cdot \av \cdot (\bv-\av)^{\rho(\pi,\ho)-1} . $$
 \label{proposition_Lambda_nu}
 \end{proposition}
 \begin{proof}
 For $\pi$ in the half open interval $[\nu,\ho)$,
 let $F(\pi)$ be the sum of the weights $\wt(x)$
 of all chains
 $x = \{\hz = \sigma_{0} < \sigma_{1} < \cdots < \sigma_{k} \}$
 such that
 $\sigma_{k} \join \nu = \pi$.
 Observe that
 $$   
     \Psi_{[\hz,\pi)} \cdot (\bv + (\av-\bv)) \cdot (\av-\bv)^{\rho(\pi,\ho)-1} 
  =
     \sum_{\nu \leq \sigma \leq \pi}
		   F(\sigma)    
   .   $$
 Note that $\bv + (\av - \bv) = \av$ is the result which records
 if the chain $x$ contains $\pi$ or not.
 Multiply by $(-1)^{\rho(\pi,\eta)}$ and sum over all
 $\pi \in [\nu,\eta]$.
 \begin{eqnarray}
 \sum_{\nu \leq \pi \leq \eta}
   \Psi_{[\hz,\pi)} \cdot \av \cdot (\av-\bv)^{\rho(\pi,\ho)-1} 
     \cdot
   (-1)^{\rho(\pi,\eta)}
   & = &
 \sum_{\nu \leq \pi \leq \eta}
    \sum_{\nu \leq \sigma \leq \pi}
	  F(\sigma)    
       \cdot
	  (-1)^{\rho(\pi,\eta)}      \nonumber \\
   & = &
 \sum_{\nu \leq \sigma \leq \eta}
	  F(\sigma)    
       \cdot
    \sum_{\sigma \leq \pi \leq \eta}
	  (-1)^{\rho(\pi,\eta)}      \nonumber \\
   & = &
 \sum_{\nu \leq \sigma \leq \eta}
	  F(\sigma)    
       \cdot
	  \delta_{\sigma,\eta}         \nonumber \\
   & = &
	  F(\eta)         .
 \label{equation_F}
 \end{eqnarray}
 Now the $\ab$-index of $\Lambda_{\nu}$ is given by the sum
 \begin{eqnarray*}
 \Psi_{\Lambda_{\nu}}
   & = &
 \sum_{\nu \leq \eta < \ho}
   F(\eta) \\
   & = &
 \sum_{\nu \leq \eta < \ho}
 \sum_{\nu \leq \pi \leq \eta}
   \Psi_{[\hz,\pi)} \cdot \av \cdot (\av-\bv)^{\rho(\pi,\ho)-1} 
     \cdot
   (-1)^{\rho(\pi,\eta)}    \\
   & = &
 \sum_{\nu \leq \pi < \ho}
   \Psi_{[\hz,\pi)} \cdot \av \cdot (\av-\bv)^{\rho(\pi,\ho)-1} 
     \cdot
 \left(
 \sum_{\pi \leq \eta < \ho}
   (-1)^{\rho(\pi,\eta)}
 \right)    \\
   & = &
 \sum_{\nu \leq \pi < \ho}
   \Psi_{[\hz,\pi)} \cdot \av \cdot (\av-\bv)^{\rho(\pi,\ho)-1} 
     \cdot
   (-1)^{\rho(\pi,\ho) - 1}    ,
 \end{eqnarray*}
 which proves the proposition.
 \end{proof}

 \begin{proposition}
 The sum of the weights of all chains $x$ in $\Lambda_{\nu}^{\prime}$
 that contain the element $*$
 is given by
 $$ 
 \sum_{\nu \leq \pi < \ho}
   \Psi_{[\hz,\pi)} \cdot 
 \left(
 (\av - \bv)^{\rho(\pi,\ho)-1}
 -
 \av \cdot (\av-\bv)^{\rho(\pi,\ho)-2} 
     \cdot
 \left(1 + (-1)^{\rho(\pi,\ho)}\right)  
 \right) 
    \cdot
  \bv  .
 $$
 \label{proposition_Lambda_nu_star}
 \end{proposition}
 \begin{proof}
 Assume that $\Lambda_{\nu}^{\prime}$
 has rank $n$
 and
 let $P$ be the poset $\Lambda_{\nu}^{\prime}$ with the coatoms
 removed. 
 The poset $P$ can be thought of as
 the rank selected poset
 $\Lambda_{\nu,\{1,\ldots,n-2\}}^{\prime}$.
 Observe that the rank of $P$ is one less than the rank
 of $\Lambda_{\nu}^{\prime}$.
 However we let $\rho$ denote the rank function of
 $\Lambda_{\nu}^{\prime}$ and not that of $P$.

 For $\pi$ in the poset $P$ such that $\nu \leq \pi < \ho$,
 let
 $G(\pi)$ be the sum of the weights over all chains
 $x = \{\hz = \sigma_{0} < \sigma_{1} < \cdots < \sigma_{k} \}$
 in the poset $P$ such that 
 $\nu \not\leq \sigma_{k}$ and $\sigma_{k} \join \nu = \pi$.
 Similarly, let
 $H(\pi)$ be the sum of the weights over all chains $x$
 in the poset $P$ such that $\sigma_{k} \join \nu = \pi$.
 Observe that $H(\pi) \cdot (\av-\bv) = F(\pi)$,
 where $F(\pi)$ is the function in
 the proof of Proposition~\ref{proposition_Lambda_nu}.
 Hence by equation~(\ref{equation_F}) we have that
 $$    H(\sigma)
   =
 \sum_{\nu \leq \tau \leq \sigma}
   \Psi_{[\hz,\tau)} \cdot \av \cdot (\av-\bv)^{\rho(\tau,\ho)-2} 
     \cdot
   (-1)^{\rho(\tau,\sigma)}  .   $$

 Now by chain counting we have the identity
 $$
       \Psi_{[\hz,\pi)} \cdot (\av - \bv)^{\rho(\pi,\ho)-1}
   =
       G(\pi)
    +
       \sum_{\nu \leq \sigma < \pi} H(\sigma)  .
 $$
 Bring $G(\pi)$ to the other side of the equation and sum
 over all $\pi$ such that $\nu \leq \pi < \ho$.
 We then obtain
 \begin{eqnarray*}
   &   &
 \sum_{\nu \leq \pi < \ho}
       G(\pi)
 -
 \sum_{\nu \leq \pi < \ho}
       \Psi_{[\hz,\pi)} \cdot (\av - \bv)^{\rho(\pi,\ho)-1} \\
   & = &
 -
 \sum_{\nu \leq \sigma < \pi < \ho}
       H(\sigma)    \\
   & = &
 -
 \sum_{\nu \leq \tau \leq \sigma < \pi < \ho}
   \Psi_{[\hz,\tau)} \cdot \av \cdot (\av-\bv)^{\rho(\tau,\ho)-2} 
     \cdot
   (-1)^{\rho(\tau,\sigma)}    \\
   & = &
 \sum_{\nu \leq \tau < \pi < \ho}
   \Psi_{[\hz,\tau)} \cdot \av \cdot (\av-\bv)^{\rho(\tau,\ho)-2} 
     \cdot
   (-1)^{\rho(\tau,\pi)}    \\
   & = &
 -
 \sum_{\nu \leq \tau < \ho}
   \Psi_{[\hz,\tau)} \cdot \av \cdot (\av-\bv)^{\rho(\tau,\ho)-2} 
     \cdot
 \left(1 + (-1)^{\rho(\tau,\ho)}\right) .
 \end{eqnarray*}
 Making a change of variable in the last sum and rearranging the sums gives:
 \begin{eqnarray*}
   &   &
 \sum_{\nu \leq \pi < \ho}
 G(\pi) \\
   & = &
 \sum_{\nu \leq \pi < \ho}
   \Psi_{[\hz,\pi)} \cdot 
 \left(
 (\av - \bv)^{\rho(\pi,\ho)-1}
 -
 \av \cdot (\av-\bv)^{\rho(\pi,\ho)-2} 
     \cdot
 \left(1 + (-1)^{\rho(\pi,\ho)}\right)  
 \right) .
 \end{eqnarray*}
 The result now follows by multiplying on the right by $\bv$.
 \end{proof}

 Following~\cite{BE}
 define $\cd$-polynomials $\alpha_n$
 by $\alpha_0 = -1$ and otherwise by
 \begin{eqnarray*}
 \alpha_{2k}
    & = &
  - \frac{1}{2} \left[ (\cv^2 - 2\dv)^{k} + 
	 \cv \cdot (\cv^2 - 2\dv)^{k-1} \cdot \cv \right], \\
 \alpha_{2k+1}
    & = &
     \frac{1}{2} \left[ (\cv^2 - 2\dv)^{k} \cdot \cv +
	 \cv \cdot (\cv^2 - 2\dv)^{k} \right],
 \end{eqnarray*}

 \begin{lemma}
 The $\cd$-polynomial $\alpha_{k}$ is given by
 $$
      \alpha_{k}
    =
      \av \cdot (\bv-\av)^{k-1} 
  +
 \left(
 (\av - \bv)^{k-1}
 -
 \av \cdot (\av-\bv)^{k-2}
     \cdot
 \left(1 + (-1)^{k}\right)
 \right) 
    \cdot
  \bv  .  
 $$
 \label{lemma_alpha}
 \end{lemma}
 \begin{proof}
 When $k$ is even we have
 \begin{eqnarray*}
   &   &
    -  \av \cdot (\av-\bv)^{k-1} 
    +  (\av-\bv)^{k-1} \cdot \bv
    -  2 \cdot \av \cdot (\av-\bv)^{k-2} \cdot \bv    \\
   & = &
    -  \av \cdot (\av-\bv)^{k-2} \cdot (\av - \bv)
    +  (\av - \bv) \cdot (\av-\bv)^{k-2} \cdot \bv
    -  2 \cdot \av \cdot (\av-\bv)^{k-2} \cdot \bv    \\
   & = &
    -  \av \cdot (\av-\bv)^{k-2} \cdot \av
    -  \bv \cdot (\av-\bv)^{k-2} \cdot \bv  \\
   & = &
    - \frac{1}{2} \cdot
 \left(
       (\av-\bv) \cdot (\av-\bv)^{k-2} \cdot (\av-\bv)
    +  (\av+\bv) \cdot (\av-\bv)^{k-2} \cdot (\av+\bv)
 \right)   \\
   & = &
    \alpha_{k} .
 \end{eqnarray*}
 For the odd case, begin by observing that
 $\av \cdot (\av-\bv)^{k-1}
   -
  \bv \cdot (\av-\bv)^{k-1}
 =
  (\av-\bv)^{k-1} \cdot \av
   -
  (\av-\bv)^{k-1} \cdot \bv$.
 By rearranging the terms we have
 $$
  \av \cdot (\av-\bv)^{k-1}
   +
  (\av-\bv)^{k-1} \cdot \bv
 =
  (\av-\bv)^{k-1} \cdot \av
   +
  \bv \cdot (\av-\bv)^{k-1}  .  $$
 Since the two sides are equal, they are also equal to their
 mean value. Thus we have
 $$
  \av \cdot (\av-\bv)^{k-1}
   +
  (\av-\bv)^{k-1} \cdot \bv
    =
 \frac{1}{2}
  \left(
  (\av+\bv) \cdot (\av-\bv)^{k-1}
   +
  (\av-\bv)^{k-1} \cdot (\av+\bv)
 \right)
    =
	 \alpha_{k} .
 $$
 \end{proof}

 We can now give an explicit formula
 for the $\cd$-index of the poset
 $\Lambda_{\nu}^{\prime}$. This formula generalizes
 Corollary~4.4 in~\cite{BE}.
 This corollary was proved using line shellings,
 thus restricting the results in~\cite{BE} to polytopes.

 \begin{theorem}
 For an Eulerian lattice $\Lambda$
 and an element $\nu$, $\hz < \nu < \ho$,
 the $\cd$-index of $\Lambda_{\nu}^{\prime}$
 is given by
 $$ 
     \Psi_{\Lambda_{\nu}^{\prime}}
   =
     \sum_{\nu \leq \pi < \ho}
		\Psi_{[\hz,\pi)}
	      \cdot
		\alpha_{\rho(\pi,\ho)}
  . $$
 \label{theorem_Lambda_nu_prime}
 \end{theorem}
 \begin{proof}
 Add
 Propositions~\ref{proposition_Lambda_nu}
 and~\ref{proposition_Lambda_nu_star}
 using Lemma~\ref{lemma_alpha}
 to simplify.
 \end{proof}

 We now outline a different proof
 of Theorem~\ref{thm-main} following 
 the argument of the polytope case
 in~\cite{BE}.
 First observe that we have the following
 corollary to the 
 decomposition theorem,
 Theorem~\ref{thm-decomp}.
 \begin{corollary}
 Let $\Lambda$ be a Gorenstein* lattice
 and let $\nu$ be an element of $\Lambda$ such that
 $\hz < \nu < \ho$. Then the following
 coefficientwise inequality holds
 for the $\cd$-indexes of $\Lambda_{\nu}^{\prime}$ and $\Lambda$:
 $$    \Psi_{\Lambda_{\nu}^{\prime}} \leq \Psi_{\Lambda}   .   $$
 \label{corollary_inequality}
 \end{corollary}
 \begin{proof}
 The poset $\Lambda$ is a subdivision
 of $\Lambda_{\nu}^{\prime}$
 by the map
 $\phi : \Lambda \rightarrow \Lambda_{\nu}^{\prime}$
 defined by
 $$
    \phi(\tau)
  =
    \begin{cases} \tau & \text{if $\tau \in \Lambda_{\nu}$,} \\
		  *    & \text{if $\tau \not\in \Lambda_{\nu}$.} 
    \end{cases}
 $$
 The map $\phi$ collapses everything outside
 $\Lambda_{\nu}$ to the element $*$.
 We only have to check that the inverse
 image of the element $[\hz,*]$ is a near
 Gorenstein* poset.
 But this is the content of 
 Lemma~\ref{lem-main} and
 the corollary follows from
 Theorem~\ref{thm-decomp},
 the decomposition theorem.
 \end{proof}

  Theorem~\ref{thm-main} is implied by the following
 three statements:
 (1) 
 Theorem~\ref{theorem_Lambda_nu_prime},
 (2)
 the inequality in Corollary~\ref{corollary_inequality},
 and
 (3)
 the equality~\cite[Proposition~4.6]{BE}:
 \begin{eqnarray*}
 \Pyr(\Psi_{[\tau,\pi)})  -  \alpha_{\rho(\tau,\pi)} 
   & = &
 \sum_{\tau < \sigma < \pi}
	\alpha_{\rho(\tau,\sigma)} \cdot \Pyr(\Psi_{[\sigma,\pi)}) ,
 \end{eqnarray*}
 where $[\tau,\pi]$ is an interval in an Eulerian poset. 
 The proof is verbatim to the proof of Theorem~5.1
 in~\cite{BE}.


 \section{Sheaves on posets}
 \label{section_sheaves}
 \setcounter{equation}{0}

 In this section we introduce sheaves on posets and prove
 Lemmas~\ref{lem-Gor-near-Gor} and~\ref{lem-main}.

 \begin{definition} A {\em sheaf} $F$ on a poset $\Pi$ consists of the data
 \begin{itemize}
 \item A real vector space $F_\sigma$ for each $\sigma\in\Pi$, called the
   stalk of $F$ at $\sigma$.
 \item Linear maps $\res^\sigma_\tau: F_\sigma\to F_\tau$ for each
   $\sigma>\tau$, satisfying the condition $\res^\tau_\nu \circ
   \res^\sigma_\tau = \res^\sigma_\nu$ whenever $\sigma>\tau>\nu$. These
   maps are called restriction maps.
 \end{itemize}
 A {\em map} of sheaves $F\to G$ is a collection of linear maps
 $F_\sigma\to G_\sigma$ commuting with the restriction maps.
 \end{definition}

 The main example of a sheaf is the constant sheaf $\RR_\Pi$ with all
 the
 stalks equal to $\RR$ and all the restriction maps equal to the
 identity. If $F$ is a
 sheaf on $\Pi$ and $S\subset \Pi$, we let $F|_S$  (or simply $F_S$) be
 the sheaf obtained from $F$ by setting all stalks at elements
 $\sigma\in\Pi\setmin S$ equal to zero. This makes sense for subsets
 $S$ satisfying the property that if $\tau<\sigma$ both lie in $S$,
 then the interval $[\tau,\sigma]$ also lies in $S$.

 For the remainder of this section we consider sheaves on simplicial
 complexes only. All sheaves are assumed to have finite dimensional stalks. 

 \begin{definition}
 Let $\Delta$ be a simplicial complex of rank $n$ and $F$ a sheaf on
 $\Delta$. The {\em cellular complex} $\Cd(F,\Delta)$ of $F$ is the
 complex 
 \[ 0\rar C^0\rar C^1 \rar \cdots \rar C^n \rar 0,\]
 where 
 \[ C^k
   =
    \bigoplus_{\onethingatopanother{x\in\Delta}{\rho(x) = n-k}}
	F_x.\]
 The maps $C^k\to C^{k+1}$ are defined by summing the restriction
 maps $\res^x_y$ with
 correct sign as in the simplicial chain complex of $\Delta$.
 \end{definition}

 If $S\subset\Delta$ is a subset, we write $\Cd(F,S)$ for the complex
 $\Cd(F|_S, \Delta)$. In other words, the terms in the complex
 $\Cd(F,S)$ are indexed by elements of $S$.

 Note that if $F=\RR_\Delta$ is the constant sheaf, then
 $\Cd(\RR_\Delta,\Delta)$ is the simplicial chain complex of
 $\Delta$. In particular,
 \[ H^i(\Cd(\RR_\Delta,\Delta)) = \widetilde{H}_{n-i-1}(\Delta,\RR).\]
 Recalling Definition~\ref{def-Gor}, we have that $\Delta$ is Gorenstein
 if and only if  
 \[ \dim H^i(\Cd(\RR_\Delta,[x,\ho))) = \begin{cases} 1 & \text{if $i=0$} \\
 0 & \text{otherwise.} \end{cases}\]
 for all $x\in\Delta$.

 \begin{definition} \label{def-Gor-sheaf}
  A sheaf $F$ on $\Delta$ is called {\em {\CM}} if 
 \[ H^i(\Cd(F,[x,\ho))) = 0\]
 for all $x\in\Delta$ and $i>0$. $F$ is called {\em Gorenstein} if, moreover, 
 \[  \dim H^0(\Cd(F,[x,\ho))) = \dim F_x.\]
 \end{definition}

 Thus, the complex $\Delta$ is Gorenstein if and only if the constant
 sheaf on it is Gorenstein. The definition of a {\CM} complex \cite{S1}
 similarly agrees with the notion of a {\CM} sheaf. In \cite{K} {\CM}
 sheaves were called semi-Gorenstein, but {\CM} is a more appropriate
 name. 

 Let $F$ be a {\CM} sheaf on $\Delta$. For any $x<y$ we can use
 projection maps to define a map of sheaves $F_{[x,\ho)} \to
   F_{[y,\ho)}$. This induces a map of cellular complexes, hence a
     map in cohomology: 
 \begin{equation} \label{eq-proj} H^0(\Cd(F,[x,\ho))) \rar
     H^0(\Cd(F,[y,\ho))) \end{equation} 
 We want to assemble the degree zero cohomologies into a sheaf on
 $\Delta$. However, the maps (\ref{eq-proj}) go in the wrong direction
 compared to restriction maps. To fix this, we take the dual vector
 spaces and dual maps. 

 \begin{definition}
 Let $F$ be a {\CM} sheaf on $\Delta$.
 Define the sheaf $F^\vee$ on $\Delta$ with the stalks
 \[ F^\vee_x = H^0(\Cd(F,[x,\ho)))^*\]
 and with the restriction maps being the duals of (\ref{eq-proj}).
 \end{definition}

 It is shown in \cite{K} that the assignment $F\mapsto F^\vee$ defines
 a contravariant functor from the category of {\CM} sheaves to the same
 category. In particular, $F^\vee$ is again {\CM}. Moreover,
 $F^{\vee\vee} \isom F$ canonically.

 Consider a short-exact sequence of sheaves on $\Delta$:
 \[ 0\rar F_1\rar F_2\rar F_3\rar 0.\]
 {}From the long-exact sequences in cohomology it follows that if either
 $F_1$ and $F_2$ or $F_1$ and $F_3$ are {\CM}, then so is the third sheaf
 and we get the short-exact dual sequence 
 \[ 0 \lar F_1^\vee \lar F_2^\vee \lar F_3^\vee \lar 0.\]
 If both $F_2$ and $F_3$ are {\CM}, then $F_1$ is {\CM} if and only if the
 dual map $F_3^\vee\to F_2^\vee$ is injective. 

 Suppose that in the exact sequence above both $F_2$ and $F_3$ are {\CM}
 and $F_1$ is supported on a subcomplex $\partial\Delta\subset\Delta$
 of rank $n-1$. Then $F_1$ is {\CM} as a sheaf on $\partial\Delta$ and we
 have the dual short exact sequence 
 \begin{equation} \label{eq-dual} 0\lar F_2^\vee \lar F_3^\vee \lar
   F_1^\vee \lar 0.\end{equation} 
 Here $F_1^\vee$ is computed by considering $F_1$ as a sheaf on
 $\partial\Delta$, hence on $\Delta$ the stalks of the dual sheaf are
 obtained from cohomology of degree $1$ rather than degree $0$. 

 We can now give a stronger characterization of Gorenstein and
 near-Gorenstein complexes on $\Delta$. 

 \begin{lemma} \label{lem-Gor-sheaf}
 A simplicial complex $\Delta$ is Gorenstein if and only if
 $\RR_\Delta$ is {\CM} and  
 \[ \RR_\Delta^\vee \isom \RR_\Delta.\]
 \end{lemma}

 \begin{proof} Note that by Definition~\ref{def-Gor} the complex
   $\Delta$ is Gorenstein if and only if it is {\CM} and
   $\RR_\Delta^\vee$ has one-dimensional stalks. We need to show that
   the restriction maps in the dual sheaf are all isomorphisms. Then via
   $\res^x_{\hz}$ all stalks are compatibly isomorphic to
   $\RR^\vee_{\Delta,\hz} = \RR$. 

 Decompose $\Delta = U\sqcup V$, where 
 \[ U = \{ x\in\Delta| \res^x_{\hz}: \RR_{\Delta,x}^\vee \rar
 \RR_{\Delta,\hz}^\vee \text{ is zero}\}.\] 
 Then all the restriction maps in $\RR_\Delta^\vee$ between stalks at
 points lying in different sets $U$ and $V$ are zero. Hence 
 \[ \RR_\Delta^\vee = \RR_\Delta^\vee|_U \oplus \RR_\Delta^\vee|_V. \]
 But then computing the dual $\RR_\Delta^{\vee\vee} = \RR_\Delta$, we
 see that it decomposes, which is impossible because all the restriction
 maps in $\RR_\Delta$ are the identity maps. Hence $U=\emptyset$.  
 \end{proof}

 \begin{lemma} \label{lem-nGor-sheaf}
 Let $\Delta$ be a simplicial complex of rank $n$ and
 $\partial\Delta\subset \Delta$ a subcomplex of rank $n-1$. Then $\Delta$ is
 near-Gorenstein with boundary $\partial\Delta$ if and only if
 $\RR_\Delta$ is {\CM} and  
 \[ \RR_\Delta^\vee \isom \RR_{\Int(\Delta)},\]
 where $\Int(\Delta) = \Delta\setmin\partial\Delta$.
 \end{lemma}
 \begin{proof}
 First we assume that the conditions given in the lemma are
 satisfied and show that then $\Delta$ is near-Gorenstein with the
 given boundary. To see that $\partial\Delta$ is Gorenstein of rank
 $n-1$, consider the exact sequence of sheaves 
 \begin{equation}\label{eq-sGor} 0\rar \RR_{\partial\Delta} \rar \RR_\Delta
   \rar \RR_{\Int{\Delta}} \rar 0.\end{equation} 
 Since $\RR_\Delta$ and $\RR_{\Int(\Delta)}$ are both {\CM}
 (they are dual to each other by assumption),
 and $\RR_{\partial\Delta}$ is supported on the rank
 $n-1$ subcomplex $\partial\Delta$, we get the dual sequence
 (\ref{eq-dual}) 
 \[ 0\rar \RR_{\partial\Delta}^\vee \rar
	  \RR_\Delta \rar \RR_{\Int(\Delta)} \rar 0.\]
 It follows that $\RR_{\partial\Delta}^\vee = \RR_{\partial\Delta}$ and
 $\partial\Delta$ is Gorenstein by the previous lemma. The remaining
 homology conditions of Definition~\ref{def-nGor} follow from the
 isomorphism $\RR_\Delta^\vee \isom \RR_{\Int(\Delta)}$. 

 Conversely, let us assume that $\Delta$ is near-Gorenstein with
 boundary $\partial\Delta$. Then $\RR_\Delta$ and
 $\RR_{\partial\Delta}$ are {\CM}, hence so is $\RR_{\Int(\Delta)}$ by the
 exact sequence (\ref{eq-sGor}) above. The dual sequence is  
 \[ 0\rar \RR_{\partial\Delta} \rar \RR_{\Int(\Delta)}^\vee \rar
 \RR_\Delta^\vee\rar 0.\] 
 The sheaf $\RR_\Delta^\vee$ has the same stalks as
 $\RR_{\Int(\Delta)}$, hence $\RR_{\Int(\Delta)}^\vee$ has the same
 stalks as~$\RR_\Delta$. It suffices to prove that
 $\RR_{\Int(\Delta)}^\vee \isom \RR_\Delta$. 

 As in the proof of the previous lemma, if $\RR_{\Int(\Delta)}^\vee$ is
 not isomorphic to $\RR_\Delta$, then it decomposes, and so does
 $\RR_{\Int(\Delta)}$ by duality. Since all the restriction maps in
 $\RR_{\Int(\Delta)}$ are identity maps, it follows that the poset
 $\Int(\Delta)$ decomposes into a disjoint union (of unrelated
 elements). But then any sheaf on $\Int(\Delta)$, including
 $\RR_\Delta^\vee$, also decomposes accordingly, which by duality
 implies that $\RR_\Delta$ decomposes. This is a contradiction. 
 \end{proof}

 \subsection{Complementary posets}

 Given a poset $\Pi$ containing a pair $(\Pi_1,\partial\Pi_1)$, define
 the complementary pair $(\Pi_2,\partial \Pi_2)$ by $\Pi_2= \Pi\setmin
 \Int(\Pi_1)$, and $\partial \Pi_2 = \partial\Pi_1$. 

 \begin{proposition}
 Let $\Delta$ be a Gorenstein simplicial complex of rank $n$ and
 $\Delta_1$ a near-Gorenstein subcomplex of the same rank and with
 boundary $\partial\Delta_1$. Then the complementary pair
 $(\Delta_2,\partial\Delta_2)$ is also near-Gorenstein of rank $n$. 
 \end{proposition}

 \begin{proof}
 Consider the exact sequence of sheaves 
 \[ 0\rar \RR_{\Delta_1} \rar \RR_\Delta \rar \RR_{\Int(\Delta_2)}\rar 0.\]
 Since the left two sheaves are {\CM}, so is the rightmost sheaf. The
 dual sequence is 
 \[ 0\lar \RR_{\Int(\Delta_1)} \lar \RR_\Delta \lar
 \RR_{\Int(\Delta_2)}^\vee \lar 0.\] 
 It follows that $\RR_{\Int(\Delta_2)}^\vee = \RR_{\Delta_2}$ and
 $\Delta_2$ is near-Gorenstein with boundary $\partial\Delta_2$. 
 \end{proof}

 Given an arbitrary poset $\Pi$ containing a pair
 $(\Pi_1,\partial\Pi_1)$, with the complementary pair
 $(\Pi_2,\partial\Pi_2)$, we
 want to assume that taking the order complexes of $\Pi_1$ and $\Pi_2$,
 we again get complementary pairs. For this it is sufficient that
 $\Pi_1$ is closed under going down: if $\sigma\in\Pi_1$ then
 $[\hz,\sigma]\subset \Pi_1$ and $\Int(\Pi_1)$ is closed under going
 up: if $\sigma\in \Int(\Pi_1)$ then $[\sigma,\ho)\subset
   \Int(\Pi_1)$. If $\Pi_1$ satisfies these conditions, then $\Pi_2$
   also does. When considering complementary pairs, we will always
   assume that one pair (hence also the other) satisfies these
   conditions.

 \begin{corollary}
 Let $\Pi$ be a Gorenstein* poset of rank $n$ and $\Pi_1$ a
 near-Gorenstein* sub-poset of the same rank and with boundary
 $\partial\Pi$. Then the complementary pair $(\Pi_2,\partial\Pi_2)$ is
 also near-Gorenstein* of rank $n$. \qed
 \end{corollary}

 Let us now prove Lemma~\ref{lem-Gor-near-Gor}. The poset
 $\Pi\setmin\{\pi\}$ is the complement of $[\hz,\pi]$ with boundary
 $[\hz,\pi)$. This latter pair is near-Gorenstein* by
   Lemma~\ref{lem-add-one}, hence so is $\Pi\setmin\{\pi\}$.  

 For the second statement of the lemma, suppose that
 $\Pi\setmin\{\pi\}$ is near-Gorenstein* with boundary
 $[\hz,\pi)$. Divide $\Delta=\cO(\Pi)$ into
   $\Delta_1=\cO(\Pi\setmin\{\pi\})$ and its complement
   $\Delta_2=\cO([\hz,\pi])$. Then both $\Delta_1$ and $\Delta_2$
   are near-Gorenstein, hence we have an exact sequence of {\CM} sheaves 
 \[ 0\rar \RR_{\Delta_1} \rar \RR_\Delta \rar \RR_{\Int(\Delta_2)} \rar 0\]
 with dual sequence
 \[ 0\lar \RR_{\Int(\Delta_1)} \lar \RR_\Delta^\vee \lar \RR_{\Delta_2} \lar 0.\]
 It follows that $\RR_\Delta^\vee$ has all stalks one-dimensional,
 hence $\Delta$ is Gorenstein by Definition~\ref{def-Gor}. \qed 

 \subsection{Lattices}

 Let us now prove Lemma~\ref{lem-main}. We need to show that for a
 Gorenstein* lattice $\Lambda$ and a nonzero element $\nu\in\Lambda$, the
 poset $\Lambda_1 = \Lambda\setmin[\nu,\ho)$ is near-Gorenstein*
   with boundary 
 \[ \partial \Lambda_1 = \{ \sigma\in\Lambda_1| \sigma\vee \nu <
 \ho\}.\] 

 Let $(\Lambda_2,\partial\Lambda_2)$ be the complementary pair. Then
 $\Int(\Lambda_2) = [\nu,\ho)$. We also denote the corresponding
   order complexes $\Delta = \cO(\Lambda)$, $\Delta_1 =
   \cO(\Lambda_1)$, $\Delta_2=\cO(\Lambda_2)$,
   $\partial\Delta_1=\partial\Delta_2 = \cO(\partial\Delta_1)$.  

 Consider the exact sequence. 
 \[ 0\rar \RR_{\Delta_1} \rar \RR_{\Delta} \rar \RR_{\Int(\Delta_2)} \rar 0.\]
 Suppose we know that $\RR_{\Int(\Delta_2)}$ is {\CM}, with the dual sheaf
 $\RR_{\Int(\Delta_2)}^\vee$ having the same
 stalks as $\RR_{\Delta_2}$, and that the dual map $\alpha:
 \RR_{\Int(\Delta_2)}^\vee \to \RR_\Delta^\vee$ is injective. Then all
 three sheaves are {\CM} and we get the dual exact sequence 
 \[ 0\lar \RR_{\Delta_1}^\vee \lar \RR_\Delta \stackrel{\alpha}{\lar}
 \RR_{\Int(\Delta_2)}^\vee \lar 0.\]  
 The sheaf $\RR_{\Int(\Delta_2)}^\vee$ being a subsheaf of the constant
 sheaf must be the constant sheaf~$\RR_{\Delta_2}$. It follows from
 this that both $\Delta_1$ and $\Delta_2$ are near-Gorenstein, hence
 $\Lambda_1$ and $\Lambda_2$ are near-Gorenstein*. 

 Let us take $x\in \Delta_2$. We need to show that the stalk
 $\RR_{\Int(\Delta_2),x}^\vee$ is well-defined (no higher cohomology in
 the cellular complex), one-dimensional, and that it injects
 into  $\RR_{\Delta,x}^\vee = \RR$ by the map $\alpha$.

 First consider the case when $x\in \Int(\Delta_2)$. Then the entire
 interval $[x,\ho)$ in $\Delta$ lies in
   $\Int(\Delta_2)$. The two complexes that compute the stalks of
   $\RR_{\Delta,x}^\vee$ and $\RR_{\Int(\Delta_2),x}^\vee$ are equal,
   with the map between them being the identity. 

 Second consider the case when $x\in\partial\Delta_2$. Then $x$ is a chain
 \[ x = \{\hz < \sigma_1<\sigma_2<\cdots< \sigma_k\},\]
 where $\sigma_i\vee \nu < \ho$,
 $\sigma_i\notin[\nu,\ho)$ for all $i$. Let  
 \[ \Pi = [x,\ho) \cap \Int(\Delta_2).\]
 An element of $\Pi$ is a refinement of $x$ that contains at least one
 $\tau\in \Int(\Lambda_2) = [\nu,\ho)$. Such a $\tau$ must
   necessarily satisfy $\sigma_k<\tau$. If we let
   $S=[\sigma_k\vee\nu,\ho)$, then we have 
 \[ \Pi = \bigcup_{\tau\in S} [\{\hz<\sigma_1<\cdots<\sigma_k<\tau\}, \ho).\]

 It remains to prove the cellular complex of
 $\RR_\Pi$ has one-dimensional cohomology in degree $0$ and  that the
 projection map $\RR_{[x,\ho)} \to \RR_\Pi$ induces a surjection in the
 cohomology.  

 Consider the complex of sheaves indexed by chains in $S$ where the maps
 are the projections with $\pm$ signs: 
 \begin{equation} \RR_\Pi \rar \bigoplus_{\sigma_k\vee\nu\leq\tau}
   \RR_{[\{x<\tau\},\ho)} \rar \bigoplus_{\sigma_k\vee \nu\leq\tau_1<\tau_2}
     \RR_{[\{x<\tau_1<\tau_2\},\ho)}\rar\cdots.\end{equation} 
 We claim that this gives a resolution of $\RR_\Pi$. Indeed, if
     $\pi\in\Pi$ is a chain that contains
     $\{\tau_1<\tau_2<\cdots<\tau_l\}$, where $\tau_i\in S$, then the stalk
     of this complex at $\pi$ is indexed by subsets of
     $\{\tau_1,\ldots,\tau_l\}$, hence is acyclic.  

 We can compute the cellular cohomology of $\RR_\Pi$ by a spectral
 sequence, where we first compute the cellular cohomology of each sheaf
 in the resolution. All sheaves in the complex have only degree zero
 cohomology of dimension one because  $\Delta$ is Gorenstein.  
 This implies that the cellular cohomology of $\RR_\Pi$ is computed by
 the complex 
 \begin{equation}\label{eq-cohom}
   \bigoplus_{\sigma_k\vee\nu\leq\tau} \RR \rar
   \bigoplus_{\sigma_k\vee\nu\leq\tau_1<\tau_2} \RR
   \rar\cdots.\end{equation} 
 Note that if we indexed the terms in this complex by chains
   $\tau_1<\tau_2<\cdots$, where $\tau_1$ is strictly greater than
   $\sigma_k\vee\nu$, then we would get the cellular complex of
   $\cO([\sigma_k\vee\nu,\ho))$. Since we also allow
   $\tau_1=\nu$, we get  the cellular complex of  
 \[ \cO([\nu,\ho))\times \{\hz,\ho\}\]
 but without the minimal element $\hz$. The cellular complex of the
 product is acyclic and the missing minimal element gives us that the
 the degree zero cohomology of the cellular complex of $\RR_\Pi$ has
 dimension one. 

 To see that we have surjectivity of 
 \[  H^0(\Cd(\RR_{[x,\ho)},\Delta)) \rar H^0(\Cd(\RR_\Pi,\Delta)),\]
 note that we can make the complex (\ref{eq-cohom}) exact by extending
 it from the left by  
 \[ H^0(\Cd(\RR_{[x,\ho)},\Delta)) \rar
   \bigoplus_{\sigma_k\vee\nu\leq\tau} \RR 
 \rar\cdots,\] 
 which comes from extending the complex of sheaves
 \[ \RR_{[x,\ho)} \rar \RR_\Pi \rar \bigoplus_{\sigma_k\vee\nu\leq\tau}
 \RR_{[\{x<\tau\},\ho)}\rar \cdots.\]  

 This finishes the proof of Lemma~\ref{lem-main} and Theorem~\ref{thm-main}.
 \qed

 \section{Appendix}

 The purpose of this appendix is to prove non-negativity of the $\cd$-index
 for Gorenstein* and near-Gorenstein* posets (Theorem~\ref{thm-pos-cd}
 and~\ref{thm-pos-cd-nGor}). The two theorems are proved in
 \cite{K} for complete and quasi-convex fans, while the generalization
 to Gorenstein* posets is left largely to the reader. We give the
 necessary details here.

 We start by recalling another construction of the $\cd$-index in terms
 of two operations $\cC$ and $\cD$ on {\CM} sheaves on a poset
 $\Pi$. 

 Let $\Pi$ be a rank  $n$ poset and $\cO(\Pi)$ its order
 complex. Define the map of posets:
 \begin{eqnarray*}
  \beta:& \cO(\Pi) &\rar \Pi \\
 & \{0<\sigma_1<\cdots<\sigma_k\} &\longmapsto \sigma_k 
 \end{eqnarray*}
 If $F$ is a sheaf on $\Pi$, the pullback $\beta^*(F)$ is a sheaf on
 $\cO(\Pi)$ with stalks
 \[  \beta^*(F)_x = F_{\beta(x)}\]
 and the obvious restriction maps.

 Let $F$ be a sheaf on $\Pi$ with finite dimensional stalks. We define
 the $\ab$-index of~$F$ as follows (recall the definition of $\wt(x)$
 in Section~\ref{sec-cd-index}):
 \[ \Psi_F = \sum_{x\in\cO(\Pi)} \wt(x) \dim \beta^*(F)_x. \]
 If this polynomial can be written in terms of variables $\cv=\av+\bv$
 and $\dv = \av\bv + \bv\av$, then we call the resulting
 $\cd$-polynomial the $\cd$-index of the sheaf $F$. It often happens
 that the $\ab$-index can be expressed only in the form
 \[ \Psi_F = f(\cv,\dv) + g(\cv,\dv)\av,\]
 where $f$ and $g$ are homogeneous $\cd$-polynomials of degree $n$ and
 $n-1$, respectively.

 Next we recall the proof of non-negativity of the $\cd$-index in case
 of fans. For fans one can define the cellular complex of a sheaf
 as in the case of simplicial complexes. With the cellular complex one
 gets the notion of {\CM} sheaves and duality. 

 Let $\Pi$ be the poset of an $n$-dimensional but not necessarily
 complete fan and $F$ a {\CM} sheaf on $\Pi$. We define two operations
 $\cC$ and $\cD$ that each produce another {\CM} sheaf on $(n-1)$- and
 $(n-2)$-dimensional fans, respectively.

 The sheaf $\cC(F)$ is simply the restriction of $F$ to the
 $(n-1)$-skeleton $\Pi^{\leq n-1}$ of $\Pi$. It is a {\CM}
 sheaf on $\Pi^{\leq n-1}$. The sheaf $\cD(F)$ is defined by choosing a
 surjective map 
 \[ \alpha: C(F)^\vee\rar C(F) \]  
 and setting $\cD(F)$ equal to the kernel of $\alpha$. Since
 $C(F)^\vee$ and $C(F)$ have isomorphic stalks on cones of maximal
 dimension $n-1$, it follows that $\cD(F)$ is a sheaf on the
 $(n-2)$-skeleton of $\Pi$. From the short-exact sequence of sheaves we
 get that $\cD(F)$ is a {\CM} sheaf on $\Pi^{\leq n-2}$.

 Since both $\cC$ and $\cD$ produce {\CM} sheaves, we can compose
 them. If $w(\cv,\dv)$ is a monomial in $\cv$ and $\dv$ of degree $m$, then
 applying $w(\cC,\cD)$ to the sheaf $F$, we get a sheaf on the
 $(n-m)$-skeleton of $\Pi$. In particular, if $m=n$, then
 $w(\cC,\cD)(F)$ is a sheaf on the zero cone.

 The following is the most general positivity result proved in
 \cite{K}:

 \begin{theorem} Let $\Pi$ be an $n$-dimensional fan and $F$ a {\CM}
   sheaf on~$\Pi$, such that the $\ab$-index of $F$ can be expressed in
   the form
 \[   \Psi_F = f(\cv,\dv) + g(\cv,\dv)\av.\]
 For any $\cd$-monomial  $w(\cv,\dv)$ of degree $n$, the coefficient of $w$
 in $f(\cv,\dv)$ is equal to the dimension of the stalk
 \[ w(\cC,\cD)(F)_\hz.\]
 In particular, $f(\cv,\dv)$ has non-negative integer coefficients.
 \end{theorem}

 Applying this theorem to the constant sheaf $\RR_\Pi$ on a complete or
 quasi-convex fan gives the non-negativity of the $\cd$-index.

 If, instead of sheaves on $\Pi$, we work with the dimensions of the stalks
 of such sheaves, then the operations $\cC$ and $\cD$ give a method for
 computing the coefficients of the $\cd$-index. It is clear that the
 same method then works for an arbitrary Eulerian poset. However, we
 need sheaves to show non-negativity of the $\cd$-index.

 For an arbitrary poset, the main difficulty in extending the
 operations $\cC$ and $\cD$ above lies in the definition of the
 cellular complex and hence the dual sheaf. Recall that the maps in the
 cellular complex are the restrictions with appropriate signs. In the
 case of fans one gets a compatible set of signs by orienting each cone and
 setting the signs equal to $\pm 1$ depending on the orientations. For an
 arbitrary poset, there may not exist such a compatible set of
 orientations. Thus, we have to do all computations involving cellular
 complexes on the order complex $\cO(\Pi)$.

 For the rest of this appendix, let $\Pi$ be a poset of rank $n$. $\Pi$
 need not be Gorenstein*. However, we require that for any
 $\sigma\in\Pi$ the subposet $[\hz,\sigma)$ must be Gorenstein*. The only
 case we will need is when $\Pi$ is the $n$-skeleton of a Gorenstein*
 poset of rank $\geq n$.

 \begin{definition} A sheaf $F$ on~$\Pi$ is {\em {\CM}} if $\beta^*(F)$
   is {\CM} on $\cO(\Pi)$.
 \end{definition}

 \begin{lemma} \label{lem-pull-back}  Let $F$ be a {\CM} sheaf on
   $\Pi$. Then there exists a 
   {\CM} sheaf $F^\vee$ on $\Pi$ such that 
 \[ \beta^*(F)^\vee \isom \beta^*(F^\vee).\]
 \end{lemma}

 \begin{proof} For $\sigma\in \Pi$, consider the fiber
   $\beta^{-1}(\sigma)$. This subposet of $\cO(\Pi)$ has the minimal
   element $x = \{\hz<\sigma\}$. We show that for any $y \in
   \beta^{-1}(\sigma)$, the restriction map in the dual sheaf 
 \begin{equation} \beta^*(F)^\vee_y \rar \beta^*(F)^\vee_x
   \label{eq-const}\end{equation} 
 is an isomorphism. We can then define
 \[ F^\vee_\sigma = \beta^*(F)^\vee_x\]
 and restriction maps for $\sigma>\tau$
 \[ F^\vee_\sigma = \beta^*(F)^\vee_{\{\hz<\sigma\}} \isom
 \beta^*(F)^\vee_{\{\hz<\tau<\sigma\}} \rar
 \beta^*(F)^\vee_{\{\hz<\tau\}} = F^\vee_\tau.\] 
 These maps are compatible and define the sheaf $F^\vee$. It is also
 easy to check that this $F^\vee$ satisfies the isomorphism in the
 statement of the lemma.

 To prove the isomorphism (\ref{eq-const}), note that the stalks of the
 dual sheaf are computed by the cellular complexes over the intervals
 $[y,\ho)$ and $[x,\ho)$, respectively. We write $[x,\ho)$ as a product
 of two posets
 \[ [x,\ho) = S_x \times T,\]
 where $S_x$ consists of all refinements of the chain $x$ with last
 element $\sigma$ and $T$ consists of all the chains between $\sigma$ and
 $\ho$. Similarly, write 
 \[ [y,\ho) = S_y \times T.\]
 Since the index set of each cellular complex is a product of two posets,
 we can write the complex as a double complex and compute its
 cohomology first along rows and then along columns. Thus, we have two
 double complexes and a map between them.  For a fixed $t =
 \{\hz<\tau_1<\cdots<\tau_k\} \in T$, the rows of the double complexes
 indexed by $S_x$ 
 and $S_y$, respectively, are the cellular complexes computing the dual of the
 constant sheaf with stalk $F_\tau$ on the poset
 $\cO([\hz,\sigma))$. By Lemma~\ref{lem-Gor-sheaf}, these two
   complexes have only degree zero cohomology isomorphic to  $F_\tau$
   and the map between the complexes induces an isomorphism of
   cohomologies. Thus, after taking cohomology along rows, we get two
   isomorphic complexes, indexed by the same set~$T$. Clearly, the
   cohomologies of these complexes are also isomorphic.
 \end{proof}

 \begin{remark}
 {\rm
 (1) Observe that there is no choice involved
 in the definition of $F^\vee$ given in the proof of
 Lemma~\ref{lem-pull-back}.
 In particular, $F\mapsto F^\vee$ is a contravariant functor on the
 category of {\CM} sheaves on $\Pi$. \\
 (2) It is easy to check that 
 \[ \dim F^\vee_\sigma = \sum_{\pi>\sigma} (-1)^{n-\rho(\pi)} \dim
 F_\pi.\]
 Thus, at least numerically, we can think of $F^\vee$ as computed by a
 cellular complex on $\Pi$, similarly to the case of fans.
 }
 \end{remark}

 Let us now turn to the construction of the operations $\cC$ and
 $\cD$. For $F$ a sheaf on~$\Pi$, let $\cC(F)$ be its restriction
 to the $(n-1)$-skeleton of $\Pi$. 

 \begin{lemma}
 If $F$ is a {\CM} sheaf on $\Pi$, then  $\cC(F)$ is
 {\CM} on $\Pi^{\leq n-1}$.
 \end{lemma}
 \begin{proof} On $\Pi$ we have an exact sequence of sheaves
 \[ 0\rar \cC(F) \rar F \rar \bigoplus_{\rho(\sigma)=n} F_\sigma \rar
 0,\]
 where $F_\sigma$ is the constant sheaf with stalk $F_\sigma$ supported
 on a single element $\sigma$. The  pullback of this sequence to $\cO(\Pi)$
 is also exact. Now $\beta^*(F)$ is {\CM} on $\cO(\Pi)$ by
 assumption. Each sheaf $\beta^*(F_\sigma)$ is a constant sheaf on
 $[\{\hz<\sigma\}, \ho) \subset \cO(\Pi)$. Since this subposet of
   $\cO(\Pi)$ is Gorenstein, the constant sheaf is {\CM}. From the
   exact sequence we get that $\cC(F)$ is {\CM} on $\Pi^{\leq n-1}$.
 \end{proof}

 To define the sheaf $\cD(F)$, we find a surjective map 
 \[ \alpha: \cC(F)^\vee \rar \cC(F) \]
 and set $\cD(F)$ equal to its kernel. Since $\alpha$ is surjective,
 $\cD(F)$ is a {\CM} sheaf on the $(n-2)$-skeleton of $\Pi$.

 The construction of $\alpha$ proceeds as in the case of fans in
 \cite{K}. Let $\sigma\in\Pi$ be an element of rank $n$ and $f\in
 F_\sigma$ a section. Then $f$ defines a map of sheaves
 \begin{eqnarray*} 
 \phi_f: \RR_{[\hz,\sigma]} &\rar F \\
 1_\sigma &\longmapsto f_\sigma.
 \end{eqnarray*}
 Restricting to the $(n-1)$-skeleton of $\Pi$ we get a map
 \[ \phi_f: \RR_{[\hz,\sigma)} \rar \cC(F).\]
 We claim that $\RR_{[\hz,\sigma)}$ is a {\CM} sheaf on $\Pi^{\leq n-1}$
 with an isomorphism 
 \[ \RR_{[\hz,\sigma)}^\vee \stackrel{\isom}{\rar} \RR_{[\hz,\sigma)}.\]
 Indeed, the pullback of $\RR_{[\hz,\sigma)}$ to $\cO(\Pi)$ is the
     constant sheaf on 
 \[ \Link_{\cO(\Pi)}(\{\hz<\sigma\}) \isom \cO([\hz,\sigma)),\]
 and the poset $[\hz,\sigma)$ is Gorenstein*.

 Since the duality is functorial, we get for each $f\in F_\sigma$ a map
 \[ \alpha_f: \cC(F)^\vee \stackrel{\phi_f^\vee}{\rar}
 \RR_{[\hz,\sigma)}^\vee \stackrel{\isom}{\rar} \RR_{[\hz,\sigma)}
 \stackrel{\phi_f}{\rar} \cC(F).\]
 We claim that a Zariski general linear combination of $\alpha_f$ (over
 all $\sigma$ and a basis of sections $f\in F_\sigma$) is a surjective
 map $\alpha$. For {\CM} 
 sheaves, it suffices to check surjectivity on maximal elements of the
 poset only. On an element $\tau$ of rank $n-1$, the map $\alpha$ is a
 composition 
 \[  F_\tau^* \stackrel{\phi^*}{\rar}
 \bigoplus_{\sigma,f} \RR^* \rar \bigoplus_{\sigma,f} \RR
 \stackrel{\phi}{\rar} F_\tau.\]
  Since the rightmost map $\phi$ is surjective, it is clear that for a general
  diagonal map $\oplus \RR^* \rar \oplus \RR$, the composition is also
  surjective. 

 {}From these two operations $\cC$ and $\cD$ we get as in the case of
 fans:

 \begin{theorem} Let $\Pi$ be a rank $n$ poset such that for every
   $\sigma\in\Pi$ the interval $[\hz,\sigma)$ is Gorenstein*. Let $F$
   be a {\CM} sheaf on~$\Pi$ such that its $\ab$-index can be expressed
   in the form
 \[   \Psi_F = f(\cv,\dv) + g(\cv,\dv)\av.\]
 Then for any $\cd$-monomial  $w(\cv,\dv)$ of degree $n$, the coefficient of $w$
 in $f(\cv,\dv)$ is equal to the dimension of the stalk
 \[ w(\cC,\cD)(F)_\hz.\]
 In particular, $f(\cv,\dv)$ has non-negative integer coefficients.
 \end{theorem}

 Applying this theorem to the constant sheaf $\RR_\Pi$ on a Gorenstein*
 poset $\Pi$ gives Theorem~\ref{thm-pos-cd}. To prove
 Theorem~\ref{thm-pos-cd-nGor}, we first apply the theorem to the
 constant sheaf  $\RR_\Pi$ on a near-Gorenstein*
 poset $\Pi$ to get non-negativity of $\Phi_\Pi(\cv,\dv)$; since
 $\partial\Pi$ is Gorenstein*, non-negativity of
 $\Psi_{\partial\Pi}(\cv,\dv)$ follows as before.

 \section{Acknowledgments}
 \setcounter{equation}{0}

 The first author thanks Mittag-Leffler Institute
 where part of this research was carried out. He was
 also supported by NSF Grant DMS-0200624.
 The second author was supported by
 NSERC grant RGPIN 283301.

 \end{document}